\documentclass[12pt]{article}
\usepackage{amsmath,amsfonts,amsthm,amssymb,amscd}
\usepackage[all]{xy}
 \usepackage[oldstylenums]{kpfonts}
   \usepackage[cal=boondoxo]{mathalfa} 
    \usepackage{hyperref} 

\newtheorem{thm}{Theorem}[section]

\newtheorem*{thm*}{Theorem}
\newtheorem*{corr*}{Corollary}
\newtheorem{lemma}[thm]{Lemma}
\newtheorem{prop}[thm]{Proposition}
\newtheorem*{prop*}{Proposition}

\newtheorem{corr}[thm]{Corollary}
\theoremstyle{definition}
\newtheorem{dfn}[thm]{Definition}
\newtheorem{ass}[thm]{Assumption}

\newtheorem{exmples}[thm]{Examples}

\theoremstyle{remark}
 
\newtheorem*{rmq}{\textit{Remark}}
\newtheorem{rmk}[thm]{\textit{Remark}}

\def\C{{\mathbf{C}}}

\newcommand{\Q}{{\mathbf{Q}}}

\newcommand{\R}{{\mathbf{R}}}

\def\AA{\mathcal{A}}

\def\ii{{\rm i}}

\newcommand{\comp}{\raise1pt\hbox{{$\scriptscriptstyle\circ$}}}
\def\del{\partial}

\renewcommand\setminus{-}
\def\into{\hookrightarrow}

 \def\lset{\{}  
\def\rset{\}}  
\def\set#1{\lset#1\rset} 
  
\def\mapright#1{\mathop{\vbox{\ialign{
                                ##\crcr
    ${\scriptstyle\hfil\;\;#1\;\;\hfil}$\crcr
 \noalign{\kern2pt\nointerlineskip}
    \rightarrowfill\crcr}}\;}}

\def\mapleft#1{\mathop{\vbox{\ialign{
                                ##\crcr
    ${\scriptstyle\hfil\;\;#1\;\;\hfil}$\crcr
 \noalign{\kern2pt\nointerlineskip}
    \leftarrowfill\crcr}}\;}}

\newcommand\gl[1]{\operatorname{GL}({#1})}

\newcommand\tr{\operatorname{Tr}}

\newcommand\End{\mathop{\rm End}\nolimits}

\def\half{\frac{1}{2}}
\title{On rigidity of locally symmetric spaces}
\author{Chris Peters}

\begin{document}
\maketitle
\section*{Introduction}
A   classical result due to Calabi and Vesentini \cite{cv} states that a  \emph{compact } locally  symmetric 
space  is  rigid, provided   all of its irreducible factors have  dimension at least 2.
This implies that such varieties (known to be algebraic) can be defined over a numberfield. This was first remarked
by Shimura in \cite{nodefs}. For a modern variant of the proof see \cite{spreads}.

Faltings \cite{arrigid} remarked that one can  show that the Kodaira-Spencer class 
for any "spread family"  of the given variety is zero which suffices for rigidity.  This is true  without any restriction on the
type of irreducible factors, and even for non-compact  locally  symmetric 
spaces. The proof uses    first of all Mumford's theory of toroidal compactifications \cite{toroidalcompactifications} of locally symmetric varieties 
together with the existence of  "good" extensions   of metric homogeneous vector bundles to these compactifications as shown in  \cite{propthm}.
The second ingredient is a careful analysis of the  extension of classical harmonic theory to a suitable   $L^2$ version.

I show in this note that the same techniques can be used to extend the results of Calabi and Vesentini to
the non-compact case.   This is stated as  Theorem~\ref{MainThm}.

Mumford's ideas are sketched in Sect.~\ref{sec:mumford} and in Sect.~\ref{sec:harmonic} I have explained the basic $L^2$--techniques 
used by Faltings. This is done in some detail  since the arguments in \cite{arrigid} are rather sketchy.

Thanks to  Christopher Deninger   for pointing out  to the reference \cite{arrigid}.

\section{Poincar\'e growth and good metrics}
\label{sec:mumford}

In this section I recall some concepts and results from \cite{propthm}. Let $X$ be a smooth quasi-projective complex variety
and let $\overline X$ be a "good" compactification: $\overline X$ is non-singular, projective and  $\del X:= \overline X\setminus X$
 a normal crossing divisor.  Hence, locally at a point of the boundary,  coordinates $(z_1,\dots,z_n)$ can be chosen such that the boundary is given by the
equation $z_1\cdots z_r=0$ and  $\del X$ can be covered by a collection of  polydisks $\Delta^n$ on which  $X$ cuts out
$(\Delta^*)^r \times \Delta^{n-r}$.  Let $\|\,\|_P$ be the  Poincar\'e norm  on such a product. Any smooth $p$ form,
say $\eta$ on $X$  is said to have \emph{Poincar\'e growth near the boundary}, if    for all tangent vectors
$\set{t_1,\dots,t_p}$ at a point of $\Delta^n\cap X$, one has the estimate $ | \eta(t_1,\cdots, t_p) |^2  \le   \text{Const. }   \|   t_1 \|_P \cdots \|  t_p  \|_P$.
This notion does not depend on choices.  By \cite[Prop. 1.1]{propthm}  such a form   defines a current on $X$.
Mumford calls  a smooth form $\omega$  on $X$ a \emph{good} form if $\omega$ as well as $d\omega$
have Poincar\'e growth near the boundary. 

Let $(E,h)$ be a hermitian holomorphic vector bundle on $X$. Recall the following definition:
\begin{dfn} \label{dfn:Chernconnexion} The \emph{Chern connection} for $(E,h)$ is  the unique  metric  connection $\nabla_E$ on $E$  whose $(0,1)$-part is 
the operator $\bar \del : \AA^0_X (E ) \to \AA^{0,1}_X(E)$ coming from  the complex structure on $E$.
\end{dfn}
Assume that $E= \overline E|_X$ where $\overline E$ is a holomorphic vector bundle on $\overline X$. 
\begin{dfn} The metric $h$ is good  relative to  $\overline E$, if locally near the boundary for every frame of $\overline E$ 
the following holds:
\begin{enumerate}
\item the matrix entries  $h_{ij}$   of $h$, respectively  $h^{-1}_{ij}$ of  $h^{-1}$,
with respect to the frame   grow at most logarithmically:  in   local coordinates $ z_1, \dots ,z_n$ as above, 
$|h_{ij}|, |h^{-1}_{ij}| \le \text{Const.}\cdot  (\log |z_1\cdots z_k|)^N$ for some integer $N$.
\item the entries of the connection matrix
$ \omega_h = \del h \cdot h^{-1}$ for the Chern  connection  are good forms.
\end{enumerate} 
\end{dfn}
By \cite[Prop. 1.3]{propthm} there is at most one extension $\overline E$ of $E$ such that $h$ is good relative to that extension.
Note also that the dual $E^*$ carries a natural metric and this metric is good relative $(\overline E)^*$.

If  $h$ is a good metric on  a vector bundle $E$ relative to an extension $\overline E$, then, by definition
 any  Chern form  calculated from the Chern connection
is  good and by \cite[Thm. 1.4]{propthm}, the class it represents,  is   the corresponding Chern class of $\overline E$.

\section{Relevant $L^2$ harmonic theory}
\label{sec:harmonic}

Let me continue with the set-up of the previous section. So $(E,h)$ is a hermitian holomorphic vector bundle on $X$
such that  $E$ is the restriction  to $X$ of a holomorphic vector bundle $\overline E$ on  $\overline X$ with the property  that $h$ is good
relative to $\overline E$. 
In addition, make the following, admittedly strong assumptions:
\begin{ass} \label{boundedsects}
1. $X$ carries a \emph{complete}  K\"ahler metric $h_X$ whose
$(1,1)$-form has  Poincar\'e growth near  $\del X$ (and hence its volume form has Poincar\'e growth).
\\
2. Smooth sections  of the bundle $\AA^k_{\overline X} (\overline E)$ of \emph{complex} $k$--forms  with values in $\overline E$ 
are bounded in the metric  induced from $h$ and $h_X$.
\end{ass}
Let me recall how to   introduces  metrics on the spaces $A^k(\overline E)$ of global \emph{complex} $k$-forms 
with values in $\overline E$.
On a fibre   $\AA^k_{X,x}(E)$ at $x\in X$ of the vector bundle $\AA^k_X(E)$, one  has a fiberwise metric induced by the  metrics $h$ and  $h_X$:
\begin{equation}\label{eqn:fibmetr}
h_x(\alpha\otimes s, \beta\otimes t)= h_X(\alpha,\beta) \, h(s,t),\quad \alpha,\beta\in \AA^{k}_{X,x },\, s,t\in E_x.
\end{equation}
Assumption 2  means that  for any two sections $\omega_i \in A^{k}(\overline E)$, $i=1, 2$ the function
$
\set{x\mapsto h_x(\omega_1,\omega_2)}  
$
is bounded on $X$.
Since by assumption 1,   the volume form for $h_X$ has Poincar\'e growth near $\del X$ it follows that the global inner product
\[
\langle \omega_1,\omega_2 \rangle = \int_{X}  h_x(\omega_1,\omega_2) \cdot \text{vol. form w.r. to $h_X$},
\quad \omega_1,\omega_2\in A^{k}_{\bar X}(\overline E)
\]
exists; in other words, one has an inclusion
\[
A^{k} (\overline E)\into L^2(X,A^k(E) )=\set{ \text{square integrable $E$-valued $k$ forms}}
\]
 and one can do harmonic theory for certain differential operators on these spaces. 
 The particular operators here are those that are  induced from  the Chern connection  $\nabla=\nabla_E$ (see Defn.~\ref{dfn:Chernconnexion}), namely
    \[
  \aligned
   \nabla: \AA_X^k(E) &\to  \AA^{k+1}_X(E),\quad \nabla^{0,1}=\bar\del, \\
            \alpha\otimes s &  \mapsto d\alpha\otimes s+ (-1)^k \alpha\otimes \nabla s.
            \endaligned
            \]
The operator  $\bar\del$, extends in the distributional sense to an operator
\[
\bar\del :L^2(X,A^{0,q}(E) ) \to L^2(X,A^{0,q+1}(E) )
\]
and since the metric on $X$ is complete and $\bar\del^2=0$, one can apply a result of Van Neumann (cf.~\ \cite[Sect. 12]{introhodge}) which
says that there is a formal adjoint operator $\bar\del^*: L^2(X,A^{0,q+1}(E) )\to L^2(X,A^{0,q}(E))$ in the sense of
distributions.  Moreover, the formal adjoint   of $\bar\del^*$ exists and   equals $\bar\del$.
These adjoints, viewed as operators on the bundles $\AA^{0,*}_X(E)$  coincide with the classical ones:
    \begin{lemma}  Let $*_E: \AA_X^{p,q}(E)\to \AA_X^{n-q,n-p}(E)$ be the fiber wise defined operator induced by the Hodge star-operator. 
    \\
    1) The formal adjoint  $\bar\del^*$ is  induced by
    $$
   -*_E\nabla ^{1,0} *_E: \AA^{0,q+1}_X(E)\to \AA^{0,q}_X(E).
   $$
    2) The formal adjoint of  $\nabla^{1,0}$ equals $(\nabla^{1,0})^* = -*_E \bar\del *_E$.
\end{lemma}
\proof Since $\bar\del = -(*_E \nabla ^{1,0} *_E)^*= -*_E ( \nabla ^{1,0})^* *_E$, the second assertion follows from the first.
The meaning of the first assertion is that for    $\omega_1\in  A^{0,q} (\overline E)$ and 
$\omega_2 \in A^{0,q+1} (\overline E)$ one has
\begin{equation}
\label{eqn:formaladj}
\langle \bar\del \omega_1,\omega_2 \rangle =- \langle \omega_1 , (*_E\nabla ^{1,0} *_E) \omega_2\rangle. 
\end{equation}
To show this, let me go through the classical calculation. 
First,  using the metric contraction
\[
\aligned
h_E: A^k(E) \otimes A^\ell(E)  &\to A^{k+\ell} \\
        (\alpha\otimes s, \beta\otimes t) &\mapsto h_E(s,t)\alpha\wedge \bar \beta
        \endaligned
        \]
one observes    the fundamental equaton
\begin{equation}
\label{eqn:Hstar}
h_E(\varphi_1,*_E \varphi_2)= h_x(\varphi_1,\varphi_2)\cdot \text{ vol. form }dV,\quad  x\in X,\quad \varphi_1,\varphi_2\in  A^k(E).
\end{equation} 
Next, the Chern connection being metric implies
that for the forms restricted to $X$ (denoted by the same symbols) one has
 \[
 h_E(\nabla\omega_1,*_E\omega_2) +(-1)^k h_E(\omega_1,\nabla(*_E\omega_2))= dh_E(\omega_1,*_E\omega_2), 
 \]
and hence, using \eqref{eqn:Hstar} and the relation $*_E\cdot *_E=(-1)^k$, one finds
\begin{equation}
\label{eqn:AdjCalc}
\bar\del h_E(\omega_1,*_E\omega_2)= \left[ h_x(\bar\del\omega_1, \omega_2) + h_x(\omega_1, (*_E \nabla^{1,0}  *_E) \omega_2)
\right ]\cdot dV .
\end{equation}
I claim that $\bar\del h_E(\omega_1,*_E\omega_2)$ is bounded near $\del X$ and that it integrates over $X$ to zero.
Assume this  for a moment.   Since the first term on the right   is bounded, the other is too.  Hence  after integration one  obtains  
\[
0= \langle \bar\del\omega_1 ,  \omega_2 \rangle + 
  \langle \omega_1 , (*_E \nabla^{1,0}*_E) \, \omega_2 \rangle
\]
and the result follows.

It remains to show the assertion about $\bar\del h_E(\omega_1,*_E\omega_2)$. Let $U_\delta$ be  a tubular neighborhood of $\del X$ with radius  $\delta$.
By Stokes' theorem,  
\begin{equation}\label{eqn:ErrorZero}
\int _X \bar\del h_E(\omega_1,*_E\omega_2) =\lim_{\delta\to 0}  \int_{\del U_\delta } h_E(\omega_1, *_E \omega_2)=0.
\end{equation}
The last equality  follows since by \eqref{eqn:Hstar} the integrand has Poincar\'e growth near the boundary and 
hence the   integral  tends to zero  (compare the proof of 
 \cite[Prop 1.2]{propthm}.
\endproof

The Laplacian $\Delta_E:= \bar\del\bar\del^*+\bar\del^*\bar \del$
preserves $L^2(A^{0,q}(X))$ and   the forms $\omega$ with $\Delta_E\omega=0$ are  by definition the \emph{harmonic} forms. Reasoning as in the classical situation (cf.~ \cite[Sect. 12]{introhodge}) one shows:

\begin{corr}\label{harmtools}
1. For all  $\omega \in  A^{0,q}_{\bar X}(\overline E)$ one has
\[
\langle \Delta_E \omega,\omega  \rangle=\langle \bar\del \omega ,\bar\del \omega \rangle+\langle \bar\del^*\omega ,\bar\del^*\omega \rangle.
\]
Hence, in the distributional sense, one has $\Delta_E\omega=0 \iff \bar\del \omega=0=\bar\del^*\omega$.\\
2.  There is an  orthogonal decomposition
\begin{equation}
\label{eqn:harmcompo}
L^2(X, \AA^{0,q}_X(E))= [\bar\del  A^{0,q-1}_X(E)]^{\rm cl} \oplus  [\bar\del^*  A^{0,q+1}_X(E)]^{\rm cl}
  \oplus \mathsf{H}^{0,q}_{(2)} (E), 
\end{equation}
where the symbol $^{\rm cl}$ stands for "topological closure" and
the symbol $\mathsf{H}_{(2)}$ stands for the harmonic $L^2$-forms,
i.e.  $L^2$-forms $\omega$ with $\Delta_E\omega=0 $ in the sense of distributions.
\end{corr}

To apply this, recall  that  by Dolbeault's theorem the cohomology group $ H^k(\overline X, \overline E)$ can be 
calculated as the $k$-th cohomology of the complex $\AA^{0,*}_{\overline X} \overline E)$.
 
\begin{prop}[\protect{\cite[Lemma 2]{arrigid}}]\label{prop:Faltngs} Assume that $\overline  E$ is a holomorphic vector bundle on $\overline X$ and that
 $(E=\overline E|_X,h)$  is a hermitian  bundle
on $X$ such that $h$ is good relative $\overline E$.  If assumption \ref{boundedsects} holds,
then there is natural \emph{injective}  homomorphism
\[
j^*_{L^2}: H^k(\overline X,\overline E)=H^k(\AA^{0,*}_{\overline X} (\overline E)) \to  \mathsf{H}^{0,k}_{(2)} (X,E) ,
\]
with target the space of $E$-valued harmonic square integrable $(0,k)$--forms.  
\end{prop}
\proof 
The map $j^*_{L^2}$ is induced from
orthogonal projection to $\mathsf{H}^k_{L^2} (E)$. The procedure is as follows.
Pick $\alpha\in  \AA^{0,k}_{\overline X}( \overline E)$ for which
$\bar\del\alpha=0$ representing  a given cohomology class $[\alpha]\in H^k(\overline X,\overline E)$. 
By assumption \ref{boundedsects},    $\beta= \alpha|X$   is  an $E$- valued $L^2$-form whose orthogonal projection
to the harmonic forms  is  $j^*_{L^2}\alpha$. One needs to verify independence of choices:
since $\bar\del\alpha=0$, one has $\bar\del\beta=0$ in the sense of currents and so, another representative for $\alpha$ leads to
a form which differs from  $\beta$ by a current of the form $\bar\del \gamma$. Hence the harmonic projection is independent of choices.

To see that it is injective, suppose that the harmonic part of $\beta$ vanishes. By \eqref{eqn:formaladj} one has 
$\langle \beta,\bar\del^*\varphi  \rangle= \langle\bar\del\beta ,\varphi\rangle=0$ and hence  $\beta$ belongs to the first summand of
\eqref{eqn:harmcompo} so that
\[
\beta =\lim_{j\to \infty}  \bar\del \gamma_j,\quad \gamma_j\in \AA^{0,k-1}_{\overline X}(\overline E).
\]
To test that this gives the zero class in $H^k(\overline X,\overline E)$, one uses  the Serre duality pairing:
\[
H^k(\overline X,\overline E) \otimes  H^{n-k}(\overline X,\Omega^n_{\overline X}\ \otimes \overline E^* )  \to H^{n,n}(\overline X)=\C
\]
as induced  by the pairing
\[
\AA^{0,k}_{\overline X} (\overline E) \otimes
           \AA^{0,n-k}_{\overline X} (\Omega^n_{\overline X}\otimes \overline E^*)  \to \AA^{n,n}_{\overline X}.
\]
To this end, consider for a   closed $\beta'\in \AA^{0,n-k}_{\overline X} (\Omega^n_{\overline X}\otimes \overline E^*) $. I claim that near $\del X$  it is bounded in norm.
To see this  let $s \in \Gamma(\overline X, \Omega^n_{\overline X} (\overline E^*))$, then, with 
$f$   a local equation for $\del X$, the product  $f \cdot s$ is  a section in the unique extension $\Omega^n(\overline X)(\log \del X)\otimes \overline E^*$ on $\overline X$ of the bundle $\Omega^n_X\otimes E^*$ on $X$ for which $h=h_X\otimes h_{E^*}$ is good. That this is the case will be shown later (Examples~\ref{ex:MainExs}.1). 
In particular, since $h( f\cdot s,  f\cdot s)= |f|^2 h(s,s)$ has  logarithmic growth near $\del X$ it follows that $h(s,s)$ and  hence also $h(\beta',\beta')$ must vanish  near $\del X$.
Hence $\beta'\in L^2(A^{0,n-k}_{X} (\Omega^n_{X}\otimes  E^*) )$.
The Serre pairing therefore is given by 
\[
(\beta, \beta'):=  \lim_{j\to \infty}  \int_{X} 
      \bar\del \gamma_j \wedge \beta'=  \lim_{j\to \infty} \lim_{\delta\to 0}  \int_{\del U_\delta } \gamma_j\wedge \beta',
      \]
      where $U_\delta$ is a tubular neighborhood of $\del X$ whose radius is $\delta$ (the last equation follows from   Stokes' theorem). Since $\beta'$ tends to zero near $\del X$, this integral vanishes.   Consequently, the cohomology class of $\beta$ is zero by Serre duality.
\endproof

I want to finish this section by showing that the Nakano inequality \cite{nakano} still holds for $E$-values harmonic $(0,q)$-forms on $X$.
To explain this, one needs some more notation. The Lefschetz operator $L$ -  which is wedging with the fundamental $(1,1)$--form for the metric $h_X$ -
preserves $L^2$--forms since the fundamental form has Poincar\'e growth near $\del X$.  Moreover, since $L$ is real,
\[
h_x(L\alpha,  \beta) dV= h_E(L\alpha, \beta)= L\alpha \wedge \overline{*\beta}= \alpha \wedge \overline {* (*^{-1}  L * \beta) }
\]
and so $\Lambda = *^{-1} L *$ is the formal adjoint of $L$. Since $*$ is an isometry, one concludes that also
$\Lambda $ preserves the $L^2$--forms.

\begin{lemma}[Nakano Inequality \cite{nakano}] \label{Nakano} Let $\omega\in  \mathsf{H}^{0,k}_{(2)} (X,E) $. With $F_h$ the curvature of the metric connection on $(E,h)$
and $\Lambda$ the formal adjoint of the  Lefschetz operator, one 
has the inequality
\[
\ii \langle \Lambda F_h \omega,\omega\rangle \ge 0.
\]
\end{lemma}
\proof For simplicity, write $\nabla^{1,0}=\del_E$ with adjoint   $\del_E^*$ . One has the K\"ahler identity (see e.g. \cite[Sect. 13]{introhodge})
\[
\aligned
\Lambda \bar\del -\bar\del\Lambda&= -\ii \del^*_E,
\endaligned
\]
which is derived in the $L^2$-setting as in the classical setting.
 Using this relation, $\bar\del\omega=0= \bar\del^*\omega  $, as well as $F_h(\omega) = \bar\del \del \omega $, one calculates
\[
\aligned
0\le \langle \del_E\omega, \del_E  \omega \rangle = 
\langle\del_E^* \del_E \omega,\omega\rangle & = \ii \langle\Lambda\bar\del \del_E \omega  -\bar\del \Lambda \del_E \omega,\omega\rangle\\
& = \ii \langle   \Lambda F_h  \omega,\omega  \rangle  - \ii \langle    \Lambda \del_E, \bar\del^*\omega\rangle\\
&=\ii  \langle   \Lambda F_h \omega ,\omega\rangle. \qedhere
\endaligned
\]
\endproof
 
 \section{The Calabi-Vesentini method in  the $L^2$--setting}

In this section I shall indicate how the method used in \cite[Sect. 7,8]{cv} to show vanishing of the groups $H^q(T_X)$ for $X$ compact can be adapted
step by step to  the non-compact setting.

Let $(X,h)$ be a K\"ahler manifold and let $T_X$ be the holomorphic tangent bundle. 
Suppose that the assumptions~\ref{boundedsects} hold. The metric $h$ induces hermitian metrics on the bundles
$\AA_X^{p,q}= \wedge^p T^*_X\otimes \wedge^q \bar T^*_X$ of forms on $X$ of type $(p,q)$.  
The Chern connection on $T_X$ is the standard Levi-Civita connection and its curvature is a global $T_X$--valued $(1,1)$--form:
\[
F_h \in  A_X^{1,1}(\End(T_X)).
\]
Using the metric one has an identification $\bar T_X^*\simeq  T_X$ and hence $F_h$ induces an endomorphism of $T_X\otimes T_X$:
\[
F_h \in T_X^*\otimes \bar T_X^* \otimes T^*_X\otimes T_X\simeq   T^*_X\otimes T^*_X\otimes T_X\otimes T_X\simeq \End (T_X\otimes T_X).
\]
One can show, using the Bianchi identity, that the resulting endomorphism vanishes on skew-symmetric tensors and hence induces
\begin{equation}
\label{eqn:Q}
Q: S^2 T_X \to S^2 T_X,\quad R= 2 \tr (Q),
\end{equation} where the function $R$ is the scalar curvature of the metric. The  operator  $Q$ is self-adjoint and hence at each $x\in X$ it has real
eigenvalues.
Let $\lambda_x$ be the smallest eigenvalue at $x$ and suppose that 
\begin{equation}
\label{eqn:lambda}
-\infty< \lambda:= \int_{x\in X }  \lambda_x  <0, \quad\text{$\lambda_x$ smallest eigenvalue  of $Q_x$}.
\end{equation}
The operator $Q$ together with the metric $h$  induces a Hermitian form  $h_Q$ on the bundles
$\AA^{0,q}(T_X)$, $q>0$  as follows:
\[
\aligned
h_Q:(\wedge^q \bar T_X^* \otimes T_X)   \otimes ( \wedge^q \bar T_X^* \otimes T_X)  &\simeq  T_X\otimes T_X \otimes(\wedge^q \bar T_X^* \otimes \wedge^q \bar T_X^*)\\
&\mapright{Q}  T_X\otimes T_X \otimes(\wedge^q \bar T_X^* \otimes \wedge^q \bar T_X^*)\to \C, 
\endaligned
\]
where the last map is induced from the hermitian metric $h$. 
If $h$  is K\"ahler-Einstein, one has \cite[Sect. 8]{cv}:
\begin{equation}
\label{eqn:CV2}
 \ii h_x(\Lambda  F \omega, \omega) = \frac{R}{2n} \|\omega\|^2 -h_Q(\omega,\omega)
\end{equation}
On the other hand, by   \cite[Lemma 3]{cv} one has the inequality
\begin{equation}
\label{eqn:CV}
h_Q(\omega,\omega) \ge \half (q+1) \lambda_x\{\omega\|^2.
\end{equation}
 In  (loc. cit.)  it is shown that  first of all  $R<0$  implies $\lambda <0$, and hence, combining \eqref{eqn:CV} and \eqref{eqn:CV2} that
 \begin{equation}
\label{eqn:CV3}
 \ii h_x(\Lambda  F \omega, \omega) \le \left( \frac{R}{2n} -  \half (q+1) \lambda_x \right)   \|\omega\|^2.
 \end{equation}
The above function is $\le 0$ whenever $  \frac{R}{2n} -  \half (q+1) \lambda<0$ and it is identically zero  if and only if $\omega=0$.
Now contrast this with the version~\ref{Nakano} of Nakano's Lemma which holds under  the assumptions of Sect.~\ref{sec:harmonic}.    The conclusion  is:
\begin{prop} \label{Mainl2vanishing} Suppose that the assumptions~\ref{boundedsects} hold for a  quasi projective K\"ahler-Einstein  manifold $(X,h)$
and its holomorphic  tangent bundle $(T_X,h)$.  Suppose    also that $R<0$, where $R$ is the scalar curvature. 

Then for all   integers $q$  for which $q< \frac{R}{n\lambda}-1$,   one has   $\mathsf{H}^{0,q}_{(2)}(X,T_X)=0$.
\end{prop}
 
  \begin{rmk} \label{NoVects} The above proof has to be modified slightly for $q=0$.
 In that case the term $h_Q(\omega,\omega)$ in \eqref{eqn:CV2} vanishes and since $R<0$ the above argument directly shows that $\mathsf{H}^{0}_{(2)}(X,T_X)=0$.
  This implies that $\bar X$ admits no vectorfields tangent to $\del X$.  
    \end{rmk}

\section{Application to locally symmetric varieties of hermitian type}

Let $G$ be a reductive $\Q$--algebraic group of hermitian type, i.e. for $K \subset G(\R)$ maximal compact, $D=G(\R)/K$ is a bounded 
symmetric domain. Fix some  neat arithmetic subgroup $\Gamma\subset G(\Q)$ and let $X= \Gamma\backslash D$ be the 
corresponding locally symmetric manifold. It is quasi-projective and by \cite{toroidalcompactifications} admits a smooth toroidal compactification $\overline X$ with boundary a 
normal crossing divisor  $\del X$. 

Let $\rho :G \to \gl E$ be a finite dimensional complex algebraic representation with $\tilde E_\rho$ the corresponding holomorphic  
vector bundle on $D$ and $E_\rho$ the bundle it defines on $X$. Fix also a $G$--equivariant hermitian 
metric $\tilde h$ on $\tilde E_\rho$ (which exists since the isotropy group of the $G(\R)$--action on  $D$ is the compact group $K$)
and write $h$ for the induced metric on $E_\rho$. By \cite[Thm. 3.1.]{propthm}, there is a unique extension of $E_\rho$ to an algebraic vectorbundle
$\overline E_\rho$ on $\overline X$ with the property that the metric $h$ is a   called good metric  for the bundle $E_\rho$
relative to $\overline E_\rho$.

For what follows it is important to observe:
\begin{lemma} \label{KEmetric}
The metric $(1,1)$-form $\omega_{h_X}$  of a K\"ahler-Einstein metric  $h_X$ has
 Poincar\'e growth   near $\del X$.
\end{lemma}
\proof
 The K\"ahler-Einstein condition means that
 \[
\omega_{h_X} = - k\cdot \ii \del\bar\del \log(\det h_X), 
\]
for some positive real constant $k$. Up to some positive constant, the right hand side can be identified   with the
first Chern form  for the canonical line bundle $\Omega^n_X$ with respect to the metric induced by $h_X$.  Since this metric
is $G(\R)$-equivariant,  it is good in Mumford's sense and so $\omega_{h_X} $ is also good. \endproof

Clearly, if this is to be useful in applications, given a  bundle (with some $G(\R)$--equivariant hermitian metric), one needs to get hold of 
the  extension making the metric good.
\begin{exmples} \label{ex:MainExs}1. Let $ E=\Omega^p_X$. Then $\overline E=  \Omega^p_{\overline X}(\log \del X)$, the bundle of $p$-forms with
at most log-poles along $\del X$. This is not trivial. See \cite[Prop. 3.4.]{propthm} where this is shown for $p=1$.
Since $\Omega^p_{\overline X}(\log \del X)=\bigwedge^p \Omega^1_{\overline X}(\log \del X)$ this  implies
the result for all $p$. In particular, smooth sections of $\Omega^1$ are bounded near $\del X$. Indeed, if $f=0$ is a local equation
for $\del X$ and $\omega$ a smooth section of $\Omega^1_X$, then $f\cdot \omega$ is a smooth section of $\Omega^1_X(\log \del X)$.
Then $\| f\cdot \omega\|^2= \|f\|^2 \|\omega\|^2$ and since  $\|\omega\|^2 \le C( \log \|f\|)^N$, $\| f\cdot \omega\|^2$ is bounded.
A similar argument holds for smooth sections of $\Omega^p_X$ and hence for sections of $\AA^{p,q}_X$.
\\
2.  One has  $\overline T_X= T_{\overline X}(-\log \del X)$, the bundle  of holomorphic vector fields on $\overline X$ 
which are tangent to the boundary $\del X$, since this  is  the dual of  the bundle  $\Omega^1_{\overline X}(\log \del X)$.
Any smooth section  of this bundle is  bounded near the boundary: its normal component tends to zero and the Poincar\'e growth
of the metric implies (by compactness of $\del X$) that  tangential component  remains bounded. \\
3. These two remarks show that the holomorphic tangent  bundle $T_X$ satisfies assumption~\ref{boundedsects} 2.
\end{exmples}

 I can finally state the  main result:
 \begin{thm} \label{MainThm} Let $(\overline X,\del X)$ as before, e.g. $X=\Gamma\backslash D$, $D=G(\R)/K$ hermitian symmetric,
 $\Gamma$ a neat arithmetic subgroup of $G(\Q)$ and $\overline X$ a  good   toroidal compactification with boundary $\del D$.
 Let $R$ be the scalar curvature of the $G(\R)$--equivariant (K\"ahler-Einstein) metric and let $\lambda$ be as before (cf. \eqref{eqn:lambda}). 
 Set  $\gamma(D):= R/n\lambda $. This is a positive integer and 
 \[
 \mathsf{H}^{0,q}_{(2)}(X,T_X)=0, \quad \text{ for all $q$ for which  } q < \gamma(D)-1.
 \]
 If no irreducible factor of $D$ has dimension $1$, one has $\gamma(D)\ge 3$.  In particular,   the resulting pairs $(\overline X,\del X)$  are  infinitesimally rigid.
 \end{thm}
 \proof
 Since $X$ admits a K\"ahler-Einstein metic $h_X$, by Lemma~\ref{KEmetric} its fundamental $(1,1)$-form has Poincar\'e growth near the boundary. 
 So the first assumption of  \ref{boundedsects} is fulfilled.
  By example~\ref{ex:MainExs}.3   the second condition is also fulfilled.
 
 In order to apply  Prop.~\ref{Mainl2vanishing}, one  observes that the K\"ahler manifold $X$ is homogeneous and that
 therefore $\lambda=\lambda_x$, $x\in X$, a constant. Since the scalar curvature of $D$ is known to be negative, this proves the result, except that $\gamma(D) $  is an  integer $\ge 2$.
 The calculation of $\gamma(D)$  is local and has been done in \cite{borel,cv} and it implies that it is an  integer $\ge 2$. Also, it is shown there  that  $\gamma(D)\ge 3$ whenever $D$ has no irreducible factor of dimension $1$.
 For details, see \cite[Sect. 3]{cv} and \cite[Sect. 2]{borel}.    See also Remark~\ref{details} below.
 
  I apply this to infinitesimal deformations  of $(\overline X,\del X)$ as follows.   As is well known, these correspond bijectively to elements of
$ H^1(\overline X, T_{\overline X}(-\log \del X))$. See e.g. \cite[Prop. 3.4.17]{sernesi}. 

Now assume that $\alpha\in  \AA^{0,1}_{\overline X}( T_{\overline X}(-\log \del X))$  represents  a given cohomology class $[\alpha]\in H^1(\overline X,T_{\overline X}(-\log \del X))$.  By  Prop.~\ref{prop:Faltngs}, the class $\beta=\alpha|X$ is an $L^2$-harmonic form and it suffices to show that $\beta=0$ which follows from the vanishing of $ \mathsf{H}^{0,1}_{(2)}(X,T_X)$.
\endproof

\begin{rmk} \label{details} For irreducible $D$  there is a table for the values of $\gamma(D)$   in \cite{cv} and \cite{borel}. I copy their  result: 
\begin{center}
\begin{tabular}{|c|c|c|c|c|c|c|}
\hline
type & $I_{p,q}$  & $II_m, m\ge 2$ & $III_m, m\ge 1$ &$IV_m, m\ge 3$& $V$ &$VI$ \\
\hline
$\gamma(D)$ & $p+q$ & $2(m-1) $ & $m+1$ & $m$ & $12$ &  $18$\\
\hline 
$\dim_\C D$ & $pq$ & $\half m(m-1)$  & $\half m(m+1)$ & $m$& $16$& $27$\\
\hline
\end{tabular}
\end{center}
If $D=D_1\times\cdots\times D_N$ is the decomposition into irreducible factors, one has $\gamma(D)=\min_j \gamma(D_j)$. One sees from this that $\gamma(D)\ge 2$
with equality precisely 
when  $D$ contains a factor of type $I_{1,1}\simeq II_2\simeq III_1$.  One also sees that the best vanishing result is for the unit ball $I_{p,1}$ where all groups vanish.

\end{rmk}

\begin{corr} Under the assumptions of Theorem~\ref{MainThm}, the pair $(\overline X,\del X)$ has a unique model over a number field.
\end{corr}
\proof This follows using spreads. For details see \cite{spreads,nodefs}.

\begin{rmq}
The above  theorem is false for  Shimura curves (one dimensional locally homogeneous algebraic manifolds). However, the corollary is true
since all Shimura curves have models  over $\overline\Q$. A proof which is a variant of the above method
 was given  in \cite{arrigid} which motivated in fact this note.
  \end{rmq}

\end{document}